\documentclass[12pt]{amsart}
\usepackage{geometry}
\usepackage{graphicx}
\usepackage{amscd}
\usepackage{amssymb}
\usepackage{amsmath}
\usepackage[latin1]{inputenc}

\theoremstyle{plain}                                       %
\newtheorem{thm}{\quad Theorem}                            %
\newtheorem{prop}[thm]{\quad Proposition}                  %
\theoremstyle{definition}                                  %
\newtheorem{rmk}[thm]{\quad Remark}                        %
\newcommand{\N}{{\Bbb N}}
\newcommand{\K}{{\Bbb K}}

\begin{document}

\textbf{Title:} "Self-inverse Sheffer sequences and Riordan
involutions."

\textbf{Authors:} Ana Luz\'{o}n* and Manuel A. Mor\'{o}n**

\textbf{Address:} *Departamento de Matem\'{a}tica Aplicada a los
Recursos Naturales. E.T.S. I. Montes. Universidad Polit\'{e}cnica
de Madrid. 28040-Madrid, SPAIN.

**Departamento de Geometria y Topologia. Facultad de Matematicas.
Universidad Complutense de Madrid. 28040- Madrid, SPAIN.

\textbf{e-mail:} *anamaria.luzon@upm.es

        **mamoron@mat.ucm.es

\textbf{All correspondence should be sent to:}

Ana Maria Luzón Cordero

e-mail: anamaria.luzon@upm.es

Address:

Departamento de Matem\'{a}tica Aplicada a los Recursos Naturales.

E.T.S.I. Montes.

Universidad Polit\'{e}cnica de Madrid.

28040-Madrid, SPAIN

Phone number: 34 913366399

Fax number: 34 915439557

\textbf{Abbreviated title:} "Self-inverse sequences and
involutions".



\newpage

\title{Self-inverse Sheffer sequences and Riordan involutions}
\author{Ana Luz\'{o}n* and Manuel A. Mor\'{o}n**}

\maketitle

\address{*Departamento de Matem\'{a}tica Aplicada a los
Recursos Naturales. E.T. Superior de Ingenieros de Montes.
Universidad Polit\'{e}cnica de Madrid. 28040-Madrid, SPAIN.}

\email{anamaria.luzon@upm.es}

\address{ **Departamento de Geometria y
Topologia. Facultad de Matematicas. Universidad Complutense de
Madrid. 28040- Madrid, SPAIN.}

 \email{mamoron@mat.ucm.es}

\vspace{1cm}

\begin{abstract}
In this short note we focus on self-inverse Sheffer sequences and
involutions in the Riordan group. We translate the results of
Brown and Kuczma on self-inverse sequences of Sheffer polynomials
to describe all involutions in the Riordan group.
\end{abstract}

\vspace{1cm}





\vspace{1cm}

Very recently, first in \cite{He-H-S} later in \cite{WangWang},
see also \cite{poly}, it has been established a very close
relation between the Sheffer group and the Riordan group, see
\cite{Sha91}, \cite{Spr94}. In fact, there is a natural
isomorphism between both groups. This means that the group
properties can be translated from one group to the other
equivalently. One of those properties is just the structure of
their finite subgroups.

In this short note we focus on self-inverse Sheffer sequences and
involutions in the Riordan group. They determine the corresponding
subgroups of order two. In fact, we translate the results of Brown
and Kuczma, \cite{Brown}, on self-inverse sequences of Sheffer
polynomials to describe all involutions in the Riordan group.
Although the translation is almost automatic we think that it is
still interesting to point it out, because of the relations to
some problems about involutions in the Riordan group posed in
\cite{Sha4} that motivated the paper \cite{Nkwanta1} and that has
been recently solved see \cite{Cheon}, \cite{CheSha},
\cite{CheSha2}. Also in pages 2264-2265 of \cite{BanPas} we have
some related results.

In some sense we want to point out that some aspects of Shapiro's
problem was solved, even before it was posed, if we reinterpret
this in terms of Sheffer sequences.

In this paper $\K$ always represents a field of characteristic
zero and $\N$ is the set of natural numbers including 0. The
notation used herein for Riordan arrays is that introduced in
\cite{teo}, see also \cite{BanPas}.

Up to the inconvenience that produce the fact that we call,
following \cite{Boas-Buck}, a generalized Appell  sequence
associated to Hadamard invertible series
$\displaystyle{h(x)=\sum_{n\geq0}h_nx^n}$ just to the sequence
obtained by multiplying by $h_n$ the $n$-term of the polynomial
sequence named by the same way in \cite{Brown}, we have the
following obvious result. The notation used below is just that
used in \cite{poly}, the operation group $\sharp_h$ is the umbral
composition as it is also described in \cite{He-H-S} for the
particular case $h(x)=e^x$, and $\star$ represents the Hadamard
product of series.

\begin{prop}
Let $N$ be  a natural number. Suppose a polynomial sequence of
Riordan type $(p_n(x))_{n\in\N}$ and
$\displaystyle{h(x)=\sum_{n\geq0}h_nx^n}$ be a series with
$h_n\neq0$ $\forall n\in\N$. Consider the Hadamard $h$-weighted
sequence $(p^h_n(x))=(p_n(x)\star h(x))$. Then, the $N$-fold
umbral composition, by means of $\sharp_h$, of the sequence
$(p^h_n(x))$ is the neutral element in $(\mathcal{R}_h, \sharp_h)$
if and only if $D^N=I$ in the Riordan group where $D=(d_{n,k})$ is
the Riordan matrix given by
$\displaystyle{p_n(x)=\sum_{k=0}^nd_{n,k}x^k}$.
\end{prop}

\begin{rmk}
Note that the sequence $e_n(x)=h_nx^n$ is the neutral element in
$(\mathcal{R}_h, \sharp_h)$.
\end{rmk}

We now translate the result in \cite{Brown} on Self-inverse
Sheffer sequences into Riordan involutions. In particular, using
the results of Section 3 in \cite{Brown} we can give a procedure
to compute all the elements $T(f\mid g)$ of the Riordan group such
that $T^2(f\mid g)=T(1\mid1)=I$. First if we impose $g=1$, then
$T^2(f\mid g)=T(1\mid1)$ if and only if $f=1$ or $f=-1$. To
construct the remaining cases we proceed in the following way:

Choose any series $\displaystyle{\phi=\sum_{n\geq0}\phi_nx^n}$
with $\phi_0\neq0$. Consider the Riordan matrix $T(1\mid\phi)$ .
As a consequence of the results in \cite{2ways} we get that
$T(1\mid A)=T^{-1}(1\mid\phi)$ where $A$ is the so called,
\cite{Spr94}, \cite{Rog78}, the $A$-sequence of $T(1\mid\phi)$. It
is clear that
$\displaystyle{\frac{x}{A}\left(\frac{x}{\phi}\right)=\frac{x}{\phi}
\left(\frac{x}{A}\right)}=x$, that is $\displaystyle{\frac{x}{A}}$
is the inverse, for the composition of the series
$\displaystyle{\frac{x}{\phi}}$. Following \cite{Brown} we take
\begin{equation}\label{E:invog}
    g=\frac{x}{\frac{x}{A}\left(-\frac{x}{\phi}\right)}
\end{equation}
Choose now any odd power series
$\displaystyle{u=\sum_{k\geq1}u_{2k-1}x^{2k-1}}$, finally take
\begin{equation}\label{E:invof}
    f=\frac{\pm x}{\frac{x}{A}\left(-\frac{x}{\phi}\right)}e^{u\left(\frac{x}{\phi}\right)}
\end{equation}
Consequently we have
\begin{prop}
Any Riordan involution different from the identity $I=T(1\mid1)$
and $-I=T(-1\mid1)$ is of the form $T(f\mid g)$ for $f$ and $g$
satisfying (\ref{E:invof}) and (\ref{E:invog}) respectively.
\end{prop}

We want to point out that the pair of series $(f,g)$ above is far
from being univocally determined by $\phi$ and $u$. For example
\begin{prop}
Let $\displaystyle{\phi=\sum_{n\geq0}\phi_nx^n}$ be a series with
$\phi_0\neq0$. Suppose that $\displaystyle{\frac{x}{A}}$ is the
compositional inverse of $\displaystyle{\frac{x}{\phi}}$
(equivalently $A$ is the $ A$-sequence of $T(1\mid \phi)$).
Suppose also that
$\displaystyle{g=\frac{x}{\frac{x}{A}\left(-\frac{x}{\phi}\right)}}$.
Then $g=-1$ if and only if $\phi(x)=\phi(-x)$ (i.e. $\phi$ is
even).
\end{prop}
\begin{proof}
Note that $\phi(x)=\phi(-x)$ if and only if the series
$\displaystyle{\frac{x}{\phi}}$ is odd.

If $\phi$ is even then $\displaystyle{\frac{x}{\phi}}$ is odd and
since $\displaystyle{\frac{x}{A}}$ is the compositional inverse of
an odd power series then $\displaystyle{\frac{x}{A}}$ is odd
itself. Consequently $g=-1$. On the other hand if $g=-1$ then
$\displaystyle{\frac{x}{A}\left(-\frac{x}{\phi}\right)=-x}$.
Composing by the right by $\displaystyle{\frac{x}{A}}$ we have
$\displaystyle{\frac{x}{A}(-x)=\frac{-x}{A(-x)}=-\frac{x}{A(x)}}$.
This implies that $\displaystyle{\frac{x}{\phi}}$ is odd and the
$\phi$ is even.
\end{proof}

In \cite{BanPas} we proved that for any $\alpha\neq0$, the Riordan
matrices $T(\pm1\mid\alpha x-1)$ are involutions, see page 2265 in
\cite{BanPas}. Now we are going to recover this result using the
construction above. In fact we will get a more general class of
Riordan involutions:

Let $\alpha\neq0$ and take
\[
\phi(x)=\frac{\alpha x}{\log(1-\alpha x)}
\]
consequently
\[
A(x)=\frac{\alpha x}{1-e^{\alpha x}}
\]
This implies that
\[
g=\frac{x}{\frac{x}{A}\left(-\frac{x}{\phi}\right)}=\alpha x-1
\]
We know that if
\[
f=\pm(\alpha x-1)e^{u\left(\frac{1}{\alpha}\log(1-\alpha
x)\right)}
\]
then $T(f\mid g)$ is an involution when $u$ is an odd series. In
particular if $u(x)=-\alpha x$ then we obtain that
$T(\pm1\mid\alpha x-1)$ is a Riordan involution.

From this point of view the fact that Pascal triangle is a
pseudo-involution, \cite{Nkwanta1}, is equivalent to the fact that
the classical Laguerre polynomials are self-inverse see
\cite{Brown}.


 {\bf Acknowledgment:} The first author was partially supported by
DGES grant MICINN-FIS2008-04921-C02-02. The second author was
partially supported by DGES grant MTM-2006-0825.


\begin{thebibliography}{99}

\normalsize



\bibitem{Boas-Buck} {R.P. Boas, R.C. Buck.} {\it Polynomial expansions of analytic
functions.} {Springer-Verlag} (1964.)



\bibitem{Brown} {J.W. Brown, M. Kuczma.} {\it Self-inverse Sheffer sequences.}
{SIAM J. Math. Anal.} {Vol. 7, No 3} {October 1976} (723-726).

\bibitem{Nkwanta1} {N. T. Cameron. and A. Nkwanta.} {\it On some (pseudo) involutions
in the Riordan Group.} {Journal of Integer Sequences.} {Vol. 8}
(2005) {Article 05.3.7}.


\bibitem{Cheon} {G.-S. Cheon, H. Kim.} {\it Simple proofs of open problems about
the structure of involutions in the Riordan group.} {Linear
Algebra Appl. }{ 428 }{ (2008) } {930-940}.

\bibitem{CheSha} {G.-S. Cheon, H. Kim, L.W. Shapiro.} {\it Riordan group involutions.} {Linear Algebra
Appl. }{ 428} { (2008) } {941-952}.

\bibitem{CheSha2} {G.-S. Cheon, S.-T. Jin, H. Kim, L.W. Shapiro.} {\it Riordan group involutions and the $\Delta$-sequence.}
{Discrete Applied Mathematics. }{ 157 (8)} { (2009) } {1696-1701
}.


\bibitem{He-H-S} {T-X. He, L.C. Hsu, P.J-S. Shiue.} {\it The Sheffer group and the Riordan group.}
{Discrete Applied Mathematics.} {155} {2007} (1895-1909).

%
%
%



\bibitem{2ways}
{A. Luzón. } {\it Iterative processes related to Riordan arrays:
The reciprocation and the inversion of power series. } { Preprint
(http://matematicas.montes.upm.es/ana/articulos/2ways.pdf)}.


\bibitem{teo}
{A. Luzón and M. A. Morón.} {\it Ultrametrics, Banach's fixed point
theorem and the Riordan group. } {Discrete Appl. Math. 156} (2008)
{2620-2635}.

\bibitem{BanPas}
{A. Luzón and M. A. Morón.} {\it Riordan matrices in the
reciprocation of quadratic polynomials. } {Linear Algebra Appl. }{
430} { (2009) } {2254-2270}.

\bibitem{poly}
{A. Luzón and M. A. Morón. } {\it Recurrence relations for
polynomial sequences via Riordan matrices. } { Preprint (ArXiv:
0904.2672) }.





\bibitem{Rog78}
{D.G. Rogers.} {\it Pascal triangles, Catalan numbers and renewal
arrays. } {Discrete Math. 22} (1978) {301-310}.








\bibitem{Sha91}
{L. W. Shapiro, S. Getu, W.J. Woan and L. Woodson.} {\it The Riordan
group. } {Discrete Appl. Math. 34} (1991) {229-239}.

\bibitem{Sha4}
{L. W. Shapiro.} {\it Some open question about random walks,
involutions, limiting distributions, and generating functions. }
{Adv. Appl. Math. 27} (2001) {585-596}.




\bibitem{Spr94}
{R. Sprugnoli.} {\it Riordan arrays and combinatorial sums. }
{Discrete Math. 132} (1994) {267-290}.


\bibitem{WangWang}
{W. Wang, T. Wang.} {\it Generalized Riordan arrays. } {Discrete
Math. 308} (2008) {6466-6500}.


\end{thebibliography}
\end{document}